\documentclass[12pt]{article}
\usepackage{graphicx}

\font\tendb=msbm10 at 12pt \font\sevendb=msbm10 at 9pt
\font\fivedb=msbm10 at 7pt
\newfam\dbfam
\textfont\dbfam=\tendb \scriptfont\dbfam=\sevendb
\scriptscriptfont\dbfam=\fivedb
\def\db{\fam\dbfam\tendb}
\font\eufm=eufm10\font\eufms=eufm10\font\eufmss=eufm10\newfam\eufam
\textfont\eufam=\eufm\scriptfont\eufam=\eufms\scriptscriptfont\eufam=\eufmss

\font\tendbb=msbm10 at 12pt \font\sevendbb=msbm7 at 9pt
\font\fivedbb=msbm5 at 6pt
\newfam\dbbfam
\textfont\dbbfam=\tendbb \scriptfont\dbbfam=\sevendbb
\scriptscriptfont\dbbfam=\fivedbb

 \def \Z {{\db Z}}
 \def \R {\hbox{\db R}}

 \def \S {S^{3}}

 \newcommand{\mod}{\mbox { mod }}

\def \fin {\hfill \framebox(7,7) }
 
\font\tenMmm=eusm10 at 12pt
\newfam\Mmmfam
\textfont\Mmmfam=\tenMmm

\def\illu #1 by #2 (#3){
  \vbox to #2{
    \hrule width #1 height 0pt depth 0pt
    \vfill
    \special{illustration #3} 
    }
  }

\begin{document}

\null \vspace{2cm}

\begin{center}
{\large {\bf Graph skein modules and symmetry
of spatial graphs}}\\
 Nafaa Chbili\\
\begin{footnotesize} Department of Mathematical Sciences\\
 College of Science\\
 UAE University\\
 E-mail: nafaachbili@uaeu.ac.ae\\
 \end{footnotesize}

\end{center}

\begin{abstract} In this paper, we  compute the  graph skein algebra of the punctured disk with two holes. Then, we apply the graph skein techniques developed here to  establish  necessary  conditions  for
a spatial graph to have a symmetry of order $p$, where $p$ is a prime.
The obstruction criteria introduced here extend some results obtained earlier for symmetric spatial graphs.
\end{abstract}
\section{Introduction}
A spatial graph in the three-sphere $\S$ is the image of an embedding of a graph $G$ in $\S$. It is worth mentioning that the graphs considered here may have multi-egdes and loops. Moreover, we assume that the valency of  each vertex  is greater than or equal to 3.  Let $p\geq 2$ be an integer, a spatial graph is
said to be $\Z_p$-symmetric if it can be isotoped into a spatial graph which is invariant by a rotation of order $p$. We consider two kinds of  $\Z_p$-symmetry. The first one is when the axis of the rotation does not intersect  the spatial graph. The second case is when the axis of the rotation  intercepts the graph only at one vertex $v$. In this paper, we address the following problem: let $G$ be a graph whose automorphism group has an element of order $p$ and let $\tilde G$ be a spatial embedding of $G$. Does  $\tilde G$ has a $\Z_p$-symmetry?\\
Przytycki \cite{Pr1}  introduced the notion of skein modules   of
three-manifolds.  Przytycki's first motivation was to extend the
definition of the polynomial invariants of links in $\S$ to links
in  other three-manifolds. Let $M$ be an oriented three-manifold. Let $\mathcal L$ be
the set of isotopy classes of framed links  in $M$. The Kauffman Bracket  skein
module of $M$ is defined to be the quotient of the free $\Z[A^{\pm 1}]$-module
generated by $\mathcal L$, by the  Kauffman bracket  skein
relations, see Section 3. \\
Inspired by the discovery of the quantum invariants of knots and links, Yamada
\cite{Ya1} introduced a polynomial invariant of  spatial graphs.
This topological invariant can be defined recursively by a family of local  relations of skein type on
planar diagrams of spatial graphs.  In spirit of the algebraic topology based on knots, we  define a
version of skein modules using embedded graphs instead of links. Here is an outline of our construction. Let $M$ be
an oriented  three-manifold.  A {\em ribbon graph} in $M$ is an oriented
surface in $M$ that retracts by deformation on a graph embedded in
$M$. Let $\mathcal G$ be the set of all isotopy classes of ribbon
graphs embedded in $M$. Let ${\mathcal R}=\Z[A^{\pm 1},
d^{-1}]$, where $d=-A^{2}-A^{-2}$, we define the {\em
graph skein module} of $M$, ${\mathcal Y}(M)$ to be  the quotient of
${\mathcal R}({\mathcal G})$ by the Yamada  relations, see
section 3.  If $M=F_{g,n}\times I$ where $F_{g,n}$ is an oriented surface of genus $g$ and having $n$ boundary components, then ${\mathcal Y}(M)$ has an algebra structure defined in the same way as in the case of the Kauffman bracket skein module \cite{Bu}. The unit of the multiplicative structure is the empty graph and the product of two ribbon graphs $G$ and $G'$ is obtained by taking a disjoint union of $G$ and $G'$ where $G$ is
pushed isotopically into $F_{g,n} \times [1/2,1]$ and $G'$ is pushed
isotopically into $F_{g,n} \times [0,1/2]$. Since the  multiplication  depends on the product structure on $M$, then
 we will denote the skein algebra of $F_{g,n}\times I$  by
$\mathcal Y$$(F_{g,n})$. In \cite{Ch1}, we proved that  the  skein algebra ${\mathcal Y}(F_{0,2})$ is isomorphic to the
polynomial algebra $\cal R$$[b]$, where $b$ is pictured in Figure 1.
\begin{center}
\includegraphics[width=4cm,height=4cm]{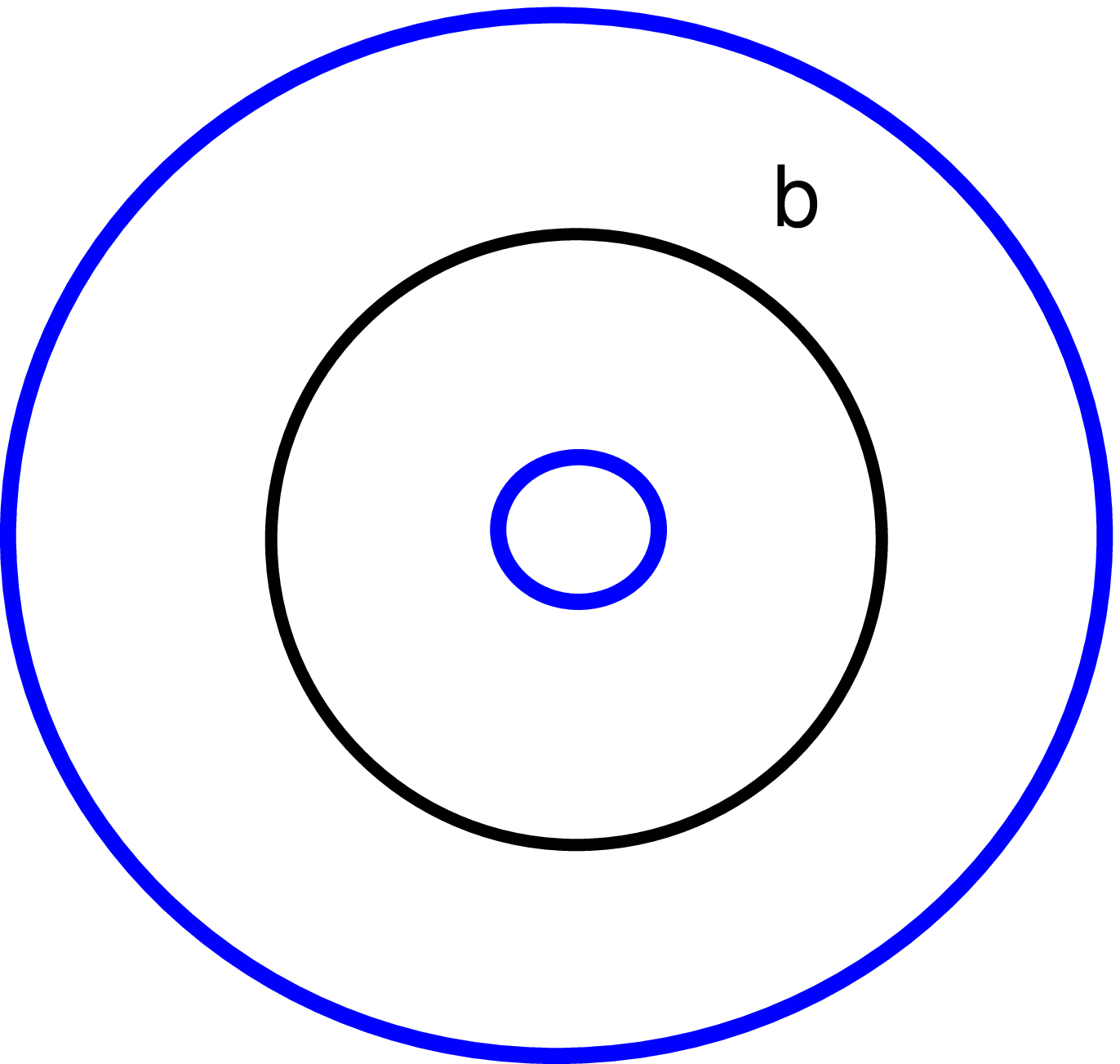}
\end{center}
\begin{center} {\sc Figure 1} \end{center}
In this paper we prove the following theorem\\
 \textbf{Theorem 1.1.}  {\sl The graph skein algebra ${\mathcal Y}(F_{0,3})$ is
isomorphic to the quotient of the  polynomial algebra $\mathcal R$$[x,y,z,t]$
by the ideal generated by
$$\begin{array}{rl}
t^2-&1+d^{-2}-2d^{-1}+(1-2d^{-1})x+(1-2d^{-1})y+z-2d^{-1}t\\
&+(1-2d^{-2})xy+xz+zy-2
d^{-1}tx-2d^{-1}ty-d^{-2}x^2-d^{-2}y^2+xyz,

\end{array}$$
where $x$, $y$, $z$ and $t$ are as in the following picture.}\\

\begin{center}
\includegraphics[width=5cm,height=4cm]{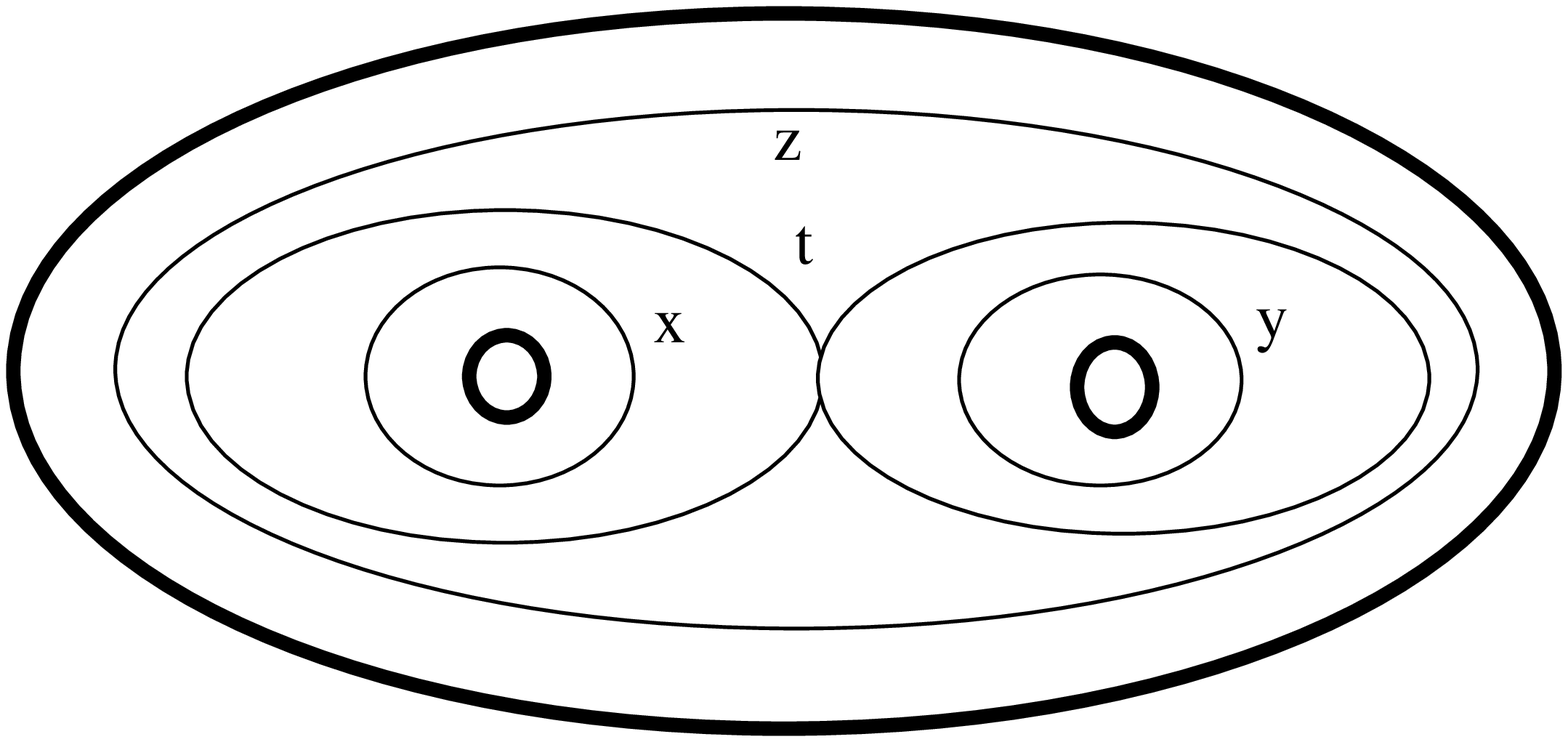}\\
Figure 2
\end{center}

Marui \cite{Ma}, studied the behavior of  the Yamada polynomial of spatial graphs with $\Z_p-$symmetry in some special cases.
 In \cite{Ch1}, we used the graph skein algebra of the annulus and the   criteria of link
periodicity introduced by Murasugi \cite{Mu}, Przytyki \cite{Pr2}
and Traczyk \cite{Tr} to  obtain an extension of Marui's
result. In the present paper, we prove some  congruence relationships satisfied by
the Yamada polynomial of spatial graphs with $\Z_p-$symmetry, where $p$ is a prime.
These relationships are, a priori,  more precise than what we obtained in  \cite{Ch1}. This is due to  the fact that  the congruence relations obtained here   involve smaller ideals than the ones in  \cite{Ch1}.
Furthermore, what we obtain here apply to vertex-fixing $\Z_p-$symmetries as well.
 This   case is not covered by the results obtained earlier in  \cite{Ch1} and \cite{Ma}. \\
\textbf{ Theorem 1.2.} {\sl Let $p$ be a prime and $\tilde G$ a
ribbon spatial   graph. Then we have the following congruences  in the ring $\Z[A^{\pm 1},d^{-1}]$\\
\texttt{ (a)} If $\tilde G$ has a $\Z_p$-symmetry with no fixed points then:  $Y(\tilde G)(A)\equiv (Y(\tilde {\underline G})(A))^p$
 modulo $p, d^{2p}-d^2$.\\
\texttt{ (b)} If $\tilde G$ has a vertex-fixing $\Z_p$-symmetry  then:  $Y(\tilde G)(A)\equiv (Y(\tilde {\underline G})(A))^p$
 modulo $p, d^{p-1}-1$.\\
\texttt{ (c)} If $\tilde G$ has a $\Z_p$-symmetry then:  $Y(\tilde G)(A)\equiv Y(\tilde G)(A^{-1})$ modulo $p, A^{8p}-1$.\\
 Where  $\tilde {\underline G}$ is the the quotient spatial graph.}\\

\textbf{Example.} Let us illustrate the last statement in  Theorem 1.2  by an example. Consider the following embedding $\tilde P$ of the Petersen graph.
 \begin{center}
\includegraphics[width=5cm,height=4cm]{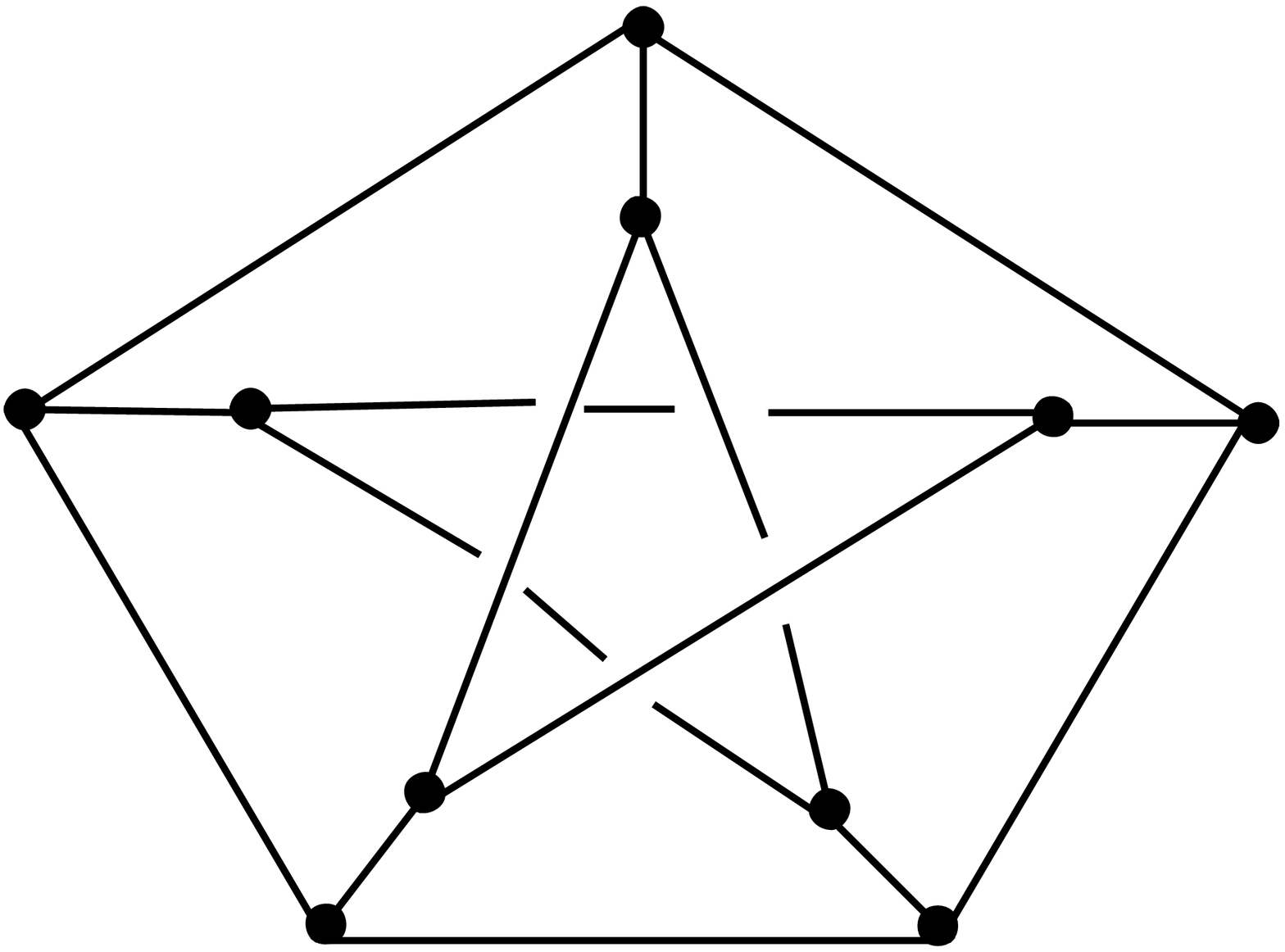}\\
\end{center}

It is known that the automorphism group of the Petersen graph is the symmetric group $S_5$. So it contains elements of order $3$ and $5$. We want  to check whether the necessary condition (c) in Theorem 1.2 can be used to rule out the possibility of symmetries of order $3$ or  $5$ for $\tilde P$.  The Yamada polynomial of $\tilde P$ is given by\\
$
\begin{array}{ll}
Y(\tilde P)(A)=& -A^{-34} - 6A^{-30} - 15A^{-26} - 35A^{-22} - 65A^{-18} - 66A^{-14} - 36A^{-10} -
15A^{-6} - 5A^{-2}\\
  & + 10 A^6 + 35 A^{10} + 61 A^{14} + 66 A^{18} + 40 A^{22} +
 15 A^{26} + 10 A^{30} + 6 A^{34} + A^{38}.

\end{array}$

Easy computation shows that  $Y(\tilde P)(A)$  is not congruent to  $Y(\tilde P)(A^{-1})$ modulo 5 and  $A^{40}-1$. Thus, $\tilde P$ doesn't have a $\Z_5$-symmetry or a  vertex fixing  $\Z_5$-symmetry. However, $Y(\tilde P)(A)$  is congruent to
 $Y(\tilde P)(A^{-1})$ modulo $A^{24}-1$ and 3. Which means that our criterion fails to decide whether $\tilde P$ has $\Z_3$-symmetry or not.
\section{The Yamada polynomial}

The  Yamada polynomial \cite{Ya1} is an invariant $R$ of regular
isotopy of spatial graphs in the three-sphere. In our terms,  $R$ is a topological invariant of ribbon spatial graphs.
It  takes its values in the ring $\Z[A^{\pm 1}]$ and may be defined recursively on planar diagrams of
spatial graphs.  Yamada also introduced a similar invariant  of trivalent graphs, with
good weight associated with the set of edges \cite{Ya2}. This construction  was extended by Yokota
\cite{Yo} using the linear skein theory introduced by Lickorich
\cite{Li}. For our reasons, we find it more convenient to
slightly change the recursive formulas introduced by Yamada, see Figure 3.
 Let $d=-A^2-A^{-2}$ and $\mathcal
R$$=\Z[A^{\pm 1}, d^{-1}]$. We define the invariant $Y$
recursively, by the four relations in the following figure. It is worth mentioning that  the following   identities hold for diagrams which are
identical except in a small disk where they look as pictured  below.\\
\null

               \begin{picture}(0,0)

               \put(135,0){\line(-1,1){13}}
               \put(119,16){\line(-1,1){13}}
               \put(105,0){\line(1,1){30}}
               \put(160,10){$= A^{4}$$Y($}
               \put(80,10){$Y($}
               \put(150,10){$)$}
               \put(205,15){\oval(28,20)[r]}
               \put(235,15){\oval(28,20)[l]}
               \put(234,10){$)$$ +A^{-4}$$Y($}
               \put(300,25){\oval(30,20)[b]}
               \put(300,0){\oval(30,20)[t]}
               \put(324,10){$)$$ -d$ $Y($}
               \put(384,0){\line(-1,1){30}}
               \put(356,0){\line(1,1){30}}
               \put(387,10){$)$}

\end{picture}
\\

\begin{picture}(0,0)
                \put(120,10){\line(-1,1){13}}
                \put(110,12){.}
                \put(110,10){.}
                \put(110,8){.}

               \put(120,10){\line(-1,-1){13}}
               \put(120,10){\line(1,0){15}}
               \put(135,10){\line(1,1){13}}
               \put(135,10){\line(1,-1){13}}
               \put(145,12){.}
                \put(145,10){.}
                \put(145,8){.}
               \put(80,10){$Y($}
               \put(150,10){$)$}
               \put(160,10){$= Y($}
                \put(200,10){\line(-1,1){13}}
               \put(200,10){\line(-1,-1){13}}
               \put(190,12){.}
                \put(190,10){.}
                \put(190,8){.}
               \put(200,10){\line(1,1){13}}
               \put(200,10){\line(1,-1){13}}
               \put(210,12){.}
                \put(210,10){.}
                \put(210,8){.}
                \put(215,10){$)-d^{-1}Y($}
                \put(275,10){\line(-1,1){13}}
               \put(275,10){\line(-1,-1){13}}
               \put(265,12){.}
                \put(265,10){.}
                \put(265,8){.}
               \put(280,10){\line(1,1){13}}
               \put(280,10){\line(1,-1){13}}
               \put(290,12){.}
                \put(290,10){.}
                \put(290,8){.}
                \put(295,10){$)$}
\end{picture}
\\

\begin{picture}(0,0)
\put(120,10){\line(-1,1){13}}
                \put(110,12){.}
                \put(110,10){.}
                \put(110,8){.}
               \put(120,10){\line(-1,-1){13}}
               \put(130,10){\circle{20}}
                \put(80,10){$Y($}
               \put(150,10){$)$}
\put(160,10){$= (d-d^{-1})Y($} \put(245,10){\line(-1,1){13}}
                \put(235,12){.}
                \put(235,10){.}
                \put(235,8){.}
               \put(245,10){\line(-1,-1){13}}
               \put(250,10){$)$}
\end{picture}
\\

\begin{picture}(0,0)

\put(80,10){$Y( \;\;D \bigsqcup \bigcirc$}
               \put(150,10){$)$}
\put(160,10){$= (d^2-1) Y(D), \mbox{ for any graph diagram } D$.}
\end{picture}
\begin{center} Figure 3
\end{center}
 The
invariant $Y$ we obtain takes values in $\mathcal
R$$=\Z[A^{\pm 1}, d^{-1}]$. It is related to Yamada's invariant $R$ by the formula:
$$Y(\tilde
G)(A)=(-d)^{\alpha(G)}R(\tilde G)(A^4),$$ where $\alpha(G)$ is equal to the number of edges of  $G$ minus the number of vertices of $G$.
\section{Graph skein modules of three-manifolds}
\subsection{The Kauffman Bracket skein module}
Let $M$ be an oriented three-manifold and $\mathcal L$ the set of
all isotopy classes of framed links in $M$. Let $\mathcal
R$=$\Z[A^{\pm 1},d^{-1}]$ and
$K(M)$ be the free $\mathcal R$-module generated by all elements
of $\mathcal L$. The Kauffman bracket skein module of
$M$ with coefficients in $\mathcal R$ which we denote here by  ${\mathcal K}(M)$ is defined as  the quotient of $K(M)$ by the
smallest submodule containing all expressions of the form:
\begin{center}
$\hspace{2cm}  \bigcirc \cup
\L - d  L $\\
$\hspace{2cm}  L -A L_{0} -A^{-1}  L_{\infty} ,$\\
\end{center}
where  $\bigcirc$ is the trivial circle and  $L$, $L_{0}$ and $ L_{\infty}$ are
three links  which are identical except in a small three-ball  where they look like  in the following picture.\\

\begin{center}
\includegraphics[width=10cm,height=2cm]{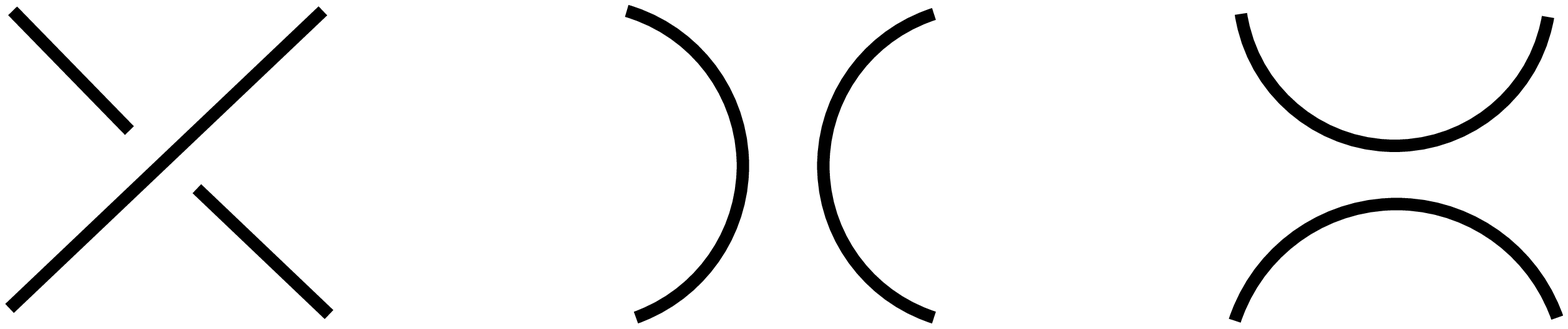} \\
{$L$}\hspace{3.5cm}{$L_0$}\hspace{3.5cm}$L_{\infty}$

\end{center}

\begin{center} {\sc Figure 4} \end{center}

The Kauffman bracket skein module has been subject to an  extensive literature. The existence and the uniqueness  of the  Kauffman bracket
polynomial \cite{Ka} is equivalent to the fact that the Kauffman
bracket skein module of $\S$ is isomorphic to $\cal R$ with the empty link $\O$ as a generator.
 The case  $M=F_{g,n}\times I$, where $F_{g,n}$ is the oriented surface of
genus $g$ with $n$ boundary components and $I$ is the unit
interval $[0,1]$ is particularly interesting. Indeed, one may project on the surface and consider diagrams of links on the surface modulo
the skein relations.  Przytycki \cite{Pr1} proved
that the Kauffman bracket skein module of $M$ is generated by all links on $F_{g,n}$ without trivial components but
including the empty link. For instance, the skein module of the solid torus $S^1\times I \times I$ is generated by $\{ b^n; n\geq 0\}$
where $b^n$ is the link in the annulus made up of $n$  parallel copies of the boundary component $b=S^1\times \{0\}\times\{0\}$, see Figure 1.\\
Furthermore, there  is an algebra   structure on the skein module of $M=F_{g,n}\times I$. The unit of the algebra is the
empty link. The product of two elements  $L$ and $L'$ is define as follows: $L.L'$ is the   disjoint union of $L$ and $L'$ where $L$ is
pushed isotopically into $F_{g,n} \times [1/2,1]$ and $L'$ is pushed
isotopically into $F_{g,n} \times [0,1/2]$. Since the  multiplication  depends on the product structure on $M$, then
 we will denote the skein algebra of $F_{g,n}\times I$  by
$\cal K$$(F_{g,n})$. Here are two examples.\\
\textbf{Theorem 3.1.1 \cite{BP}.}
{\sl 1) The algebra  ${\mathcal K}(F_{0,2})$ is isomorphic to the polynomial algebra ${\mathcal R}[b]$.\\
2) The algebra ${\mathcal K}(F_{0,3})$ is isomorphic to the polynomial algebra ${\mathcal R}[x,y,z]$, where $x, y$ and $z$ are the
three boundary components}.

\subsection{The graph skein module}
This subsection is devoted to introduce  the theory  of graph skein modules of three manifolds. This construction is inspired by  Przytycki's algebraic topology based on knots. Namely, we extend Przytycki's construction
to a similar theory  using ribbon embedded graphs instead of framed links and Yamada relations, Figure 5, instead of the Kauffman bracket skein relations. Here is an outline of the theory  and some basic results.\\
 Let $M$ be an oriented three-manifold and let $\mathcal G$ be the set of
all embeddings of  ribbon graphs in $M$ considered up to isotopy. Let ${\mathcal
R}{\mathcal G}$ be the free ${\mathcal R}$-module generated by
${\mathcal G}$. The graph skein module of $M$ with coefficients in $\mathcal R$,  ${\mathcal Y}(M,{\mathcal R},A)$ is defined to be  the quotient
 of  ${\mathcal R}{\mathcal G}$ by the smallest submodule containing all expressions of the form:\\
\null

               \begin{picture}(0,0)

               \put(135,0){\line(-1,1){13}}
               \put(119,16){\line(-1,1){13}}
               \put(105,0){\line(1,1){30}}
               \put(160,10){$- A^{4}$}
\put(205,15){\oval(28,20)[r]}
               \put(235,15){\oval(28,20)[l]}

               \put(234,10){$ -A^{-4}$}
               \put(300,25){\oval(30,20)[b]}
               \put(300,0){\oval(30,20)[t]}

               \put(324,10){$ +\;d$}
               \put(384,0){\line(-1,1){30}}
               \put(356,0){\line(1,1){30}}

\end{picture}
\\

\begin{picture}(0,0)
                \put(120,10){\line(-1,1){13}}
                \put(110,12){.}
                \put(110,10){.}
                \put(110,8){.}

               \put(120,10){\line(-1,-1){13}}
               \put(120,10){\line(1,0){15}}
               \put(135,10){\line(1,1){13}}
               \put(135,10){\line(1,-1){13}}
               \put(145,12){.}
                \put(145,10){.}
                \put(145,8){.}

               \put(160,10){$-$}
                \put(200,10){\line(-1,1){13}}
               \put(200,10){\line(-1,-1){13}}
               \put(190,12){.}
                \put(190,10){.}
                \put(190,8){.}
               \put(200,10){\line(1,1){13}}
               \put(200,10){\line(1,-1){13}}
               \put(210,12){.}
                \put(210,10){.}
                \put(210,8){.}
                \put(215,10){$\;+\;d^{-1}$}
                \put(275,10){\line(-1,1){13}}
               \put(275,10){\line(-1,-1){13}}
               \put(265,12){.}
                \put(265,10){.}
                \put(265,8){.}
               \put(280,10){\line(1,1){13}}
               \put(280,10){\line(1,-1){13}}
               \put(290,12){.}
                \put(290,10){.}
                \put(290,8){.}

\end{picture}
\\

\begin{picture}(0,0)
\put(120,10){\line(-1,1){13}}
                \put(110,12){.}
                \put(110,10){.}
                \put(110,8){.}
               \put(120,10){\line(-1,-1){13}}
               \put(130,10){\circle{20}}

\put(160,10){$-(d-d^{-1})$}

\put(245,10){\line(-1,1){13}}
                \put(235,12){.}
                \put(235,10){.}
                \put(235,8){.}
               \put(245,10){\line(-1,-1){13}}

\end{picture}
\\

\begin{picture}(0,0)

\put(120,10){$D \bigsqcup \bigcirc$}
 \put(160,10){$- (d^2-1)D , \mbox{ for any ribbon graph } D$}
\end{picture}

\begin{center} {\sc  Figure 5} \end{center}

 Throughout the rest of the paper and when there is no need to specify the value of $A$, we write ${\mathcal Y}(M)$ instead of
 ${\mathcal Y}(M,{\mathcal R},A)$. The  fact that the Yamada polynomial is uniquely determined by the Yamada skein
   relations \cite{Ya1} translates in the language of graph skein modules as follows:\\

 {\bf Theorem 3.2.1.} {\sl The graph skein module ${\mathcal Y}(S^3)$ is isomorphic to $\mathcal R$ with the empty graph $\O$ as a generator.}

In \cite{Ch1}, we have studied  the skein algebra of the annulus and proved the following:\\
{\bf Theorem 3.2.2.}  {\sl The  skein algebra ${\mathcal Y}(F_{0,2})$ is isomorphic to the
polynomial algebra $\cal R$$[b]$, where $b$ is pictured in Figure 1.}\\

Before we prove Theorem 1.1, we will  explore the  relationship between our graph-algebra and the
Kauffman bracket skein algebra. We first introduce  the Jones-Wenzel idempotents \cite{KL}.
 Let $n$ be an integer, an $n$-tangle $T$ is a one-dimensional
sub-manifold of $\R^2\times I$, such that the boundary  of $T$ is
made up of 2$n$ points $\{(0,i,1),(0,i,0)\;0\leq i\leq n-1\}$. Let
${\mathcal T}_n$ be the free $\mathcal R$-module generated by the
set of all $n$-tangles. We define $\tau _n$ to be the quotient of
${\mathcal T}_n$ by the Kauffman bracket skein relations. It is
well known that $\tau _n$ is isomorphic to the Temperley-Lieb
algebra. A set of  generators
$(U_i)_{0\leq i \leq n-1}$ of $\tau _n$ is given as follows.\\

\begin{center}

\includegraphics[width=10cm,height=4cm]{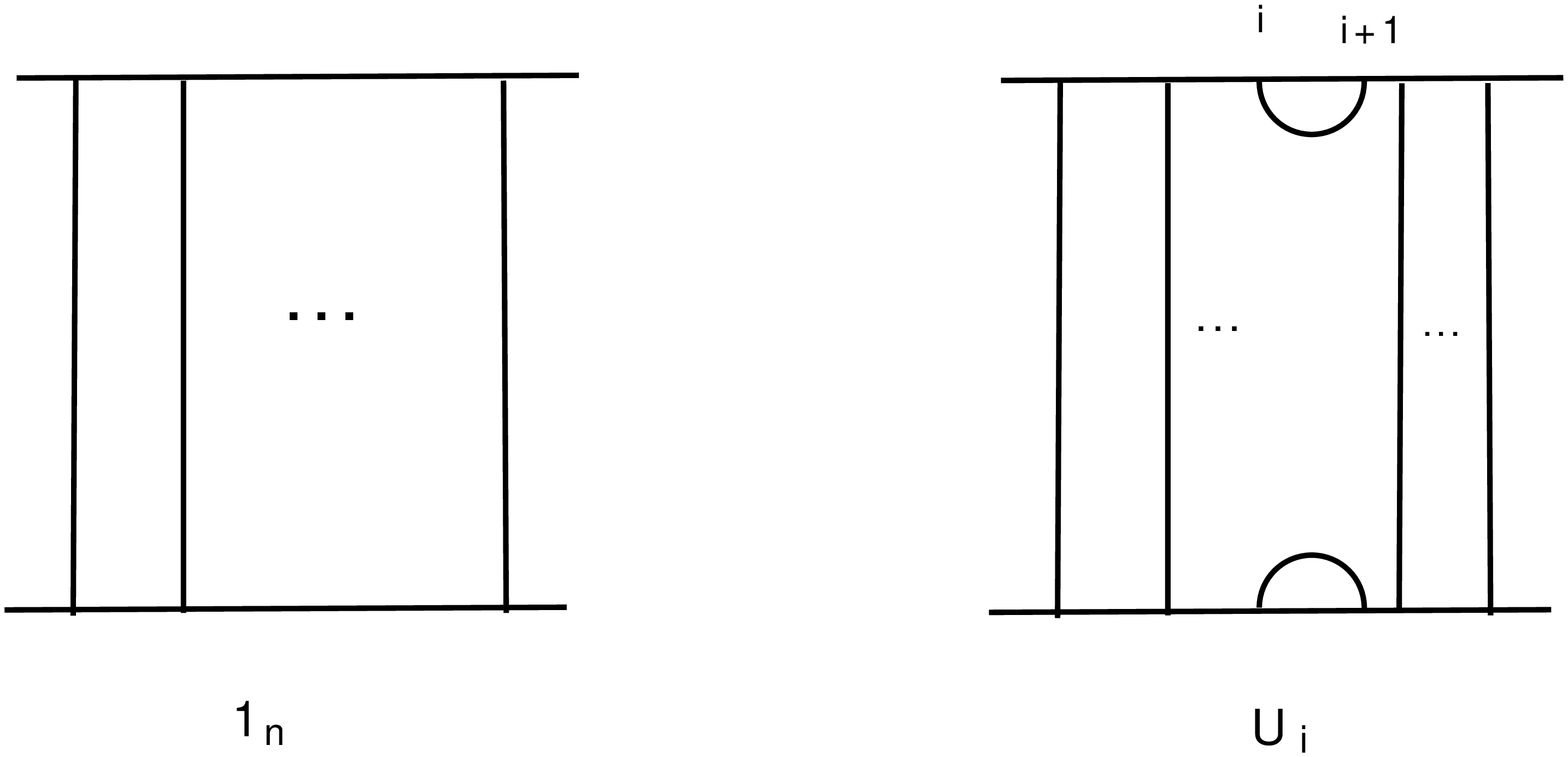}

 {\sc Figure 6} \end{center}

 Let $(f_i)_{0\leq i\leq n-1}$ denote the family of
Jones-Wenzl projectors in $\tau _n$. This family is  defined by
the
following recursive formulas:\\
$f_0=1_n,$\\
$f_{k+1}=f_k-\mu _{k+1}f_kU_{k-1}f_k$,\\
where $\mu _1=d^{-1}$ and $\mu _{k+1}=(d-\mu _k)^{-1}$.\\
In particular, we have  $f_1$=$1_n-d^{-1}U_1$. The elements
$f_k$ enjoy the following properties: $f_{k}^2=f_k$ and
$f_iU_j=U_jf_i=0$ for
$j\leq i$. See \cite{KL} for more details.\\

Let $G$ be a graph diagram lying in a given oriented surface.  Let $f_{1}$ be the Jones-Wenzl projector   in $\tau_2$
\begin{center}
\begin{picture}(0,0)
\put(0,0){$f_{1}=$} \put(30,-10){\line(0,1){25}}
\put(37,-10){\line(0,1){25}} \put(40,0){$-d^{-1}$}
\put(75,14){\oval(7,17)[b]} \put(75,-10){\oval(7,17)[t]}
\end{picture}\\
\end{center}

 We define $G'$ to be the linear combination  of graph
diagrams obtained from $G$ by replacing each edge of $G$ by two
planar strands with a projector $f_1$ in the cable, and by
replacing each vertex of $G$ by a diagram as follows\\
\begin{center}
\includegraphics[width=8cm,height=3cm]{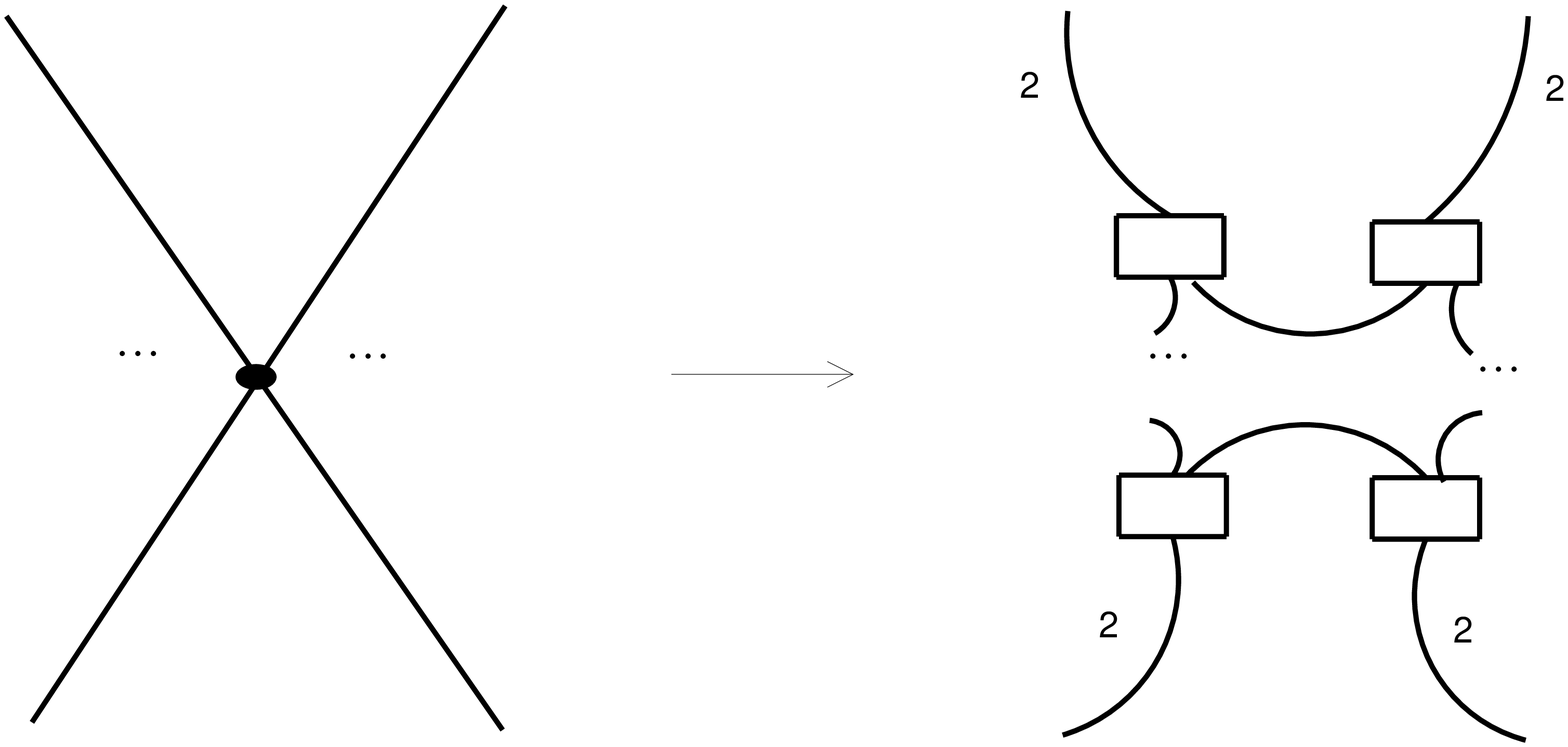}
\end{center}
\begin{center} {\sc  Figure 7} \end{center}

    Here, writing an integer $n$ beneath an edge $e$ means that this
edge has to be replaced by $n$ parallel ones. \\
 Now, let $\varphi$ be the map from ${\mathcal R}({\mathcal G})$ to $\mathcal
 K$$(F)$ defined on the generators by $\varphi (G)=G'$ and
 extended by linearity to ${\mathcal R}({\mathcal G})$, see also \cite{Ya2} and
 \cite{Yo}.\\
It was proved in \cite{Ch1} that  $\varphi$ defines a map $\Phi$ from the graph skein module  ${\mathcal Y}(F \times  I)$
to the Kauffman bracket skein module ${\mathcal K}(F \times I)$,  where $F$ is any oriented surface. Obviously,  $\Phi$ is a homomorphism of algebras. In the case of the annulus, we have seen that  both ${\mathcal K}(F_{0,2})$ and ${\cal Y}(F_{0,2})$ are isomorphic to the polynomial algebra
${\cal R}[b]$.  We can easily see  that $\Phi (b)=b^2-1$. Thus, $\Phi$ is injective. Actually, $\Phi$  defines
an isomorphism between the graph algebra ${\cal Y}(F_{0,2})$ and the
even part of the Kauffman bracket skein algebra ${\cal K}(F_{0,2})$.
\subsection{Proof of Theorem 1.1} By the reduction arguments used in the case of the annulus \cite{Ch1},  we can  prove that any graph diagram in $F_{0,3}$ can be
written as a linear combination of finite disjoint unions of bouquets each of which has  no contractible cycles. Arguments similar to the ones
used to prove  Lemma 3.3.2 in \cite{Ch1} would enable us to prove that
the graph skein module of $F_{0,3}\times I$ is generated by elements of type
$x^iy^jz^kt^{\epsilon}$, where $i,j$ and  $k$ are nonnegative
integers, and $\epsilon \in\{0,1\}$.\\
To prove that the algebra is exactly as depicted in  Theorem 1.1, we
use the fact that the homomorphism  $\Phi$ is injective. The
injectivity of $\Phi$ is due to the following identities:\\
$\Phi(x)=x^2-1$, $\Phi(y)=y^2-1$, $\Phi(z)=z^2-1$ and\\
$\Phi(t)=xyz-d^{-1}x^2-d^{-1}y^2+d^{-1}$.
Let ${\mathcal E}=<x^2,y^2,z^2, xyz>$ be the even part of the algebra  ${\mathcal K}(F_{0,3})$.
Then, $\Phi$ defines an isomorphism between ${\mathcal Y}(F_{0,3})$ and ${\mathcal E}$. The inverse isomorphism $\Psi$ is defined
as follows:
$\Psi(x^2)=x+1$, $\Psi(y^2)=y+1$, $\Psi(z^2)=z+1$ and\\
$\Psi(xyz)=t+d^{-1}x+d^{-1}y+d^{-1}$.

Using the relations above we can see that:
$$xyz=\Phi(t+d^{-1}x+d^{-1}y+d^{-1}).$$
Hence:\\
$\Phi(t^2)=\Phi(t)^2=\Phi((x+1)(y+1)(z+1)+d^{-2}(x+1)^2+d^{-2}(y+1)^2+d^{-2}$\\
 $-2d^{-1}(t+d^{-1}x+d^{-1}y+1)(x+y+1)
+2d^{-1}(x+1)(y+1)-2d^{-2}(x+y+2) ).$ \\
Since $\Phi$ is
injective, then
$$\begin{array}{rl}
t^2=&1+d^{-2}-2d^{-1}+(1-2d^{-1})x+(1-2d^{-1})y+z-2d^{-1}t\\
&+(1-2d^{-2})xy+xz+zy-2
d^{-1}tx-2d^{-1}ty-d^{-2}x^2-d^{-2}y^2+xyz.
\end{array}$$
This ends the proof of Theorem 1.1. \fin\\

In a joint work, T. Fleming and the author \cite{CF} studied the graph skein algebras of the torus $F_{1,0}$ and the punctured torus $F_{1,1}$.
They  determined a set of generators for each of these algebras.  It was proved that the graph skein algebra ${\mathcal Y}(F_{1,0})$ is
 generated by the three torus
 curves $(1,0), (0,1), (1,1)$ and the  wedge  $(1,0)\vee(0,1)$. A similar statement  was proved for the punctured torus.
\section{Proof of Theorem 1.2}
Statement (a) in Theorem 1.2  is concerned with the case of  $\Z_p-$symmetries in which the spatial graph does not intersect
 the axis of the rotation. Such a spatial graph
is called $p-$periodic. Marui \cite{Ma}, used the Yamada
polynomial to study  the periodicity of  spatial graphs with wrapping  number 1 or 2. In \cite{Ch1}, we used the criteria of link
periodicity introduced by Murasugi \cite{Mu}, Przytyki \cite{Pr2}
and Traczyk \cite{Tr} to  obtain a generalization of Marui's result. We proved the following  \cite{Ch1}: \\

\textbf{Theorem 4.1.} {\sl Let $p$  be a prime  and $\tilde G$ a
ribbon spatial graph. If $\tilde G$ is $p-$periodic, then}\\
\begin{tabular}{ll}
\texttt{(a)}& $Y(\tilde G)(A)\equiv (Y(\tilde {\underline G})(A))^p$
 modulo $p, d^{p}-d$.\\
\texttt{(b)}& $Y(\tilde G)(A)\equiv Y(\tilde G)(A^{-1})$ modulo $p,
A^{2p}-1$.
\end{tabular}\\
{\sl Where the congruences hold in the ring $\Z[A^{\pm 1},d^{-1}]$.}\\

 Now, we shall start the  proof of the first congruence relation in Theorem 1.2. We shall then explain how to extend the result to vertex-fixing $\Z_p$-symmetry as in statement (b). The idea  of the proof is to change the coefficients in the Yamada skein relations in order to define a kind of equivariant graph skein module. We already know that this idea works  well for
 the study of symmetries of links \cite{Ch2,Ch3}. Let ${\mathcal R}_p=\Z/p\Z[A^{\pm1},d^{-1}]$ and let  ${\mathcal
S}_p$ be  the free ${\mathcal R}_p$-module generated by all isotopy classes of  ribbon
graphs embedded in the solid torus. Now, let ${\mathcal Q}_p$ be
the submodule of ${\mathcal S}_p$
generated by all elements of the form:\\

\null

               \begin{picture}(0,0)

               \put(135,0){\line(-1,1){13}}
               \put(119,16){\line(-1,1){13}}
               \put(105,0){\line(1,1){30}}
               \put(160,10){$- A^{4p}$}
\put(205,15){\oval(28,20)[r]}
               \put(235,15){\oval(28,20)[l]}

               \put(234,10){$ -A^{-4p}$}
               \put(300,25){\oval(30,20)[b]}
               \put(300,0){\oval(30,20)[t]}

               \put(324,10){$ +\;d^p$}
               \put(384,0){\line(-1,1){30}}
               \put(356,0){\line(1,1){30}}

\end{picture}
\\

\begin{picture}(0,0)
                \put(120,10){\line(-1,1){13}}
                \put(110,12){.}
                \put(110,10){.}
                \put(110,8){.}

               \put(120,10){\line(-1,-1){13}}
               \put(120,10){\line(1,0){15}}
               \put(135,10){\line(1,1){13}}
               \put(135,10){\line(1,-1){13}}
               \put(145,12){.}
                \put(145,10){.}
                \put(145,8){.}

               \put(160,10){$-$}
                \put(200,10){\line(-1,1){13}}
               \put(200,10){\line(-1,-1){13}}
               \put(190,12){.}
                \put(190,10){.}
                \put(190,8){.}
               \put(200,10){\line(1,1){13}}
               \put(200,10){\line(1,-1){13}}
               \put(210,12){.}
                \put(210,10){.}
                \put(210,8){.}
                \put(215,10){$\;+\;d^{-p}$}
                \put(275,10){\line(-1,1){13}}
               \put(275,10){\line(-1,-1){13}}
               \put(265,12){.}
                \put(265,10){.}
                \put(265,8){.}
               \put(280,10){\line(1,1){13}}
               \put(280,10){\line(1,-1){13}}
               \put(290,12){.}
                \put(290,10){.}
                \put(290,8){.}

\end{picture}
\\

\begin{picture}(0,0)
\put(120,10){\line(-1,1){13}}
                \put(110,12){.}
                \put(110,10){.}
                \put(110,8){.}
               \put(120,10){\line(-1,-1){13}}
               \put(130,10){\circle{20}}

\put(160,10){$-(d-d^{-1})^p$}

\put(245,10){\line(-1,1){13}}
                \put(235,12){.}
                \put(235,10){.}
                \put(235,8){.}
               \put(245,10){\line(-1,-1){13}}

\end{picture}
\\

\begin{picture}(0,0)

\put(120,10){$D \bigsqcup \bigcirc$}
 \put(160,10){$- (d^2-1)^p D , \mbox{ for any graph diagram } D$}
\end{picture}

\begin{center} {\sc  Figure 8} \end{center}
We define ${\mathcal Y}_p$ to be the quotient of ${\mathcal S}_p$
by the submodule ${\mathcal Q}_p$. By  universal coefficient property of the skein modules, we can prove that  the module ${\mathcal Y}_p$ is isomorphic to ${\mathcal R}_p[b]$. Let $\pi:
S^1\times D^2 \longrightarrow S^1\times D^2$ denote the $p$-fold
cyclic cover defined by the action of the rotation $h$ on the solid
torus. Let $F$ (resp. $F'$) be the map from ${\mathcal S}_p$ to
${\mathcal R}_p$ defined on the set of generators of ${\mathcal
S}_p$ by $F(g)= Y(\pi ^{-1}(g))$ (resp. $F'( g)=(Y(g))^p$) and
extended to ${\mathcal S}_p$ by linearity. Using the fact that $p$ is a prime and that the finite cyclic group of order $p$
acts semi-freely on the set of states of the Yamada resolution of the diagram of the periodic spatial graph   ${\tilde G}=\pi^{-1}(g)$, we should be able
to easily prove
 the following lemma (see also  Lemma 3.4 in  \cite{Pr2}).\\
{\bf Lemma 4.2.} {\sl $F({\mathcal Q}_p)=F'({\mathcal Q}_p)=0$.}\\
According to this lemma, $F$ (resp. $F'$) defines a map  from the
skein module ${\mathcal Y}_p$ to ${\mathcal R}_p$. We will denote
this map by $\bar F$ (resp. $\bar F'$).\\
The module ${\mathcal Y}_p$ is generated by $\{b^k, k\geq 0\}$.
Let $I$ be the submodule generated by $\{\bar F(b^k)-\bar F'(b^k), k\geq 0\}$. Simple computations show that $I$ is equal
to the ideal generated by $(d^2-1)^{p}-(d^2-1)$.
Consequently: $$\bar F(g) \equiv  \bar F'(g) \mod { (d^2-1)^{p}-(d^2-1) }.$$
This ends the proof of statement (a).\\
Now we shall  consider the case of vertex fixing $\Z_p$-symmetry. Assume that $\tilde G$ is a spatial graph which is invariant by a rotation of order $p$ such that the axis of the rotation intercepts the graph only at one vertex $v$. Take a diagram of $\tilde G$ which is invariant by a planar rotation centered at $v$. If we use Yamada relations to reduce the graph diagram of $\tilde G$, then we only  need to consider equivariant states, since the contribution of non-equivariant  states  sums to zero modulo $p$. Each resolution is actually a diagram which is made up of bouquets. We distinguish two cases:\\
- If the diagram of the  state  does not contain any bouquet centered at the vertex $v$. In this case, we have a periodic diagram and we use the statement (a) to conclude that the contribution of the state is congruent to the contribution of the quotient state modulo $p$, $d^{2p}-d^2$.\\
- If the diagram of the state   contains a bouquet $B$ of $kp$-leaves centered at $v$, then the quotient state should contain a bouquet $\bar B$ with  $k$-leaves centered at $v$.
Since $Y(B)=(d-d^{-1})^{kp-1}(d^2-1)$ and  $Y(\bar B)=(d-d^{-1})^{k-1}(d^2-1)$, then the contribution of $B$ to $Y(\tilde G)$ and the contribution
of    of $\bar B$ to $Y(\underline{\tilde G})^p$ are congruent modulo $p$ and $d^{p-1}-1$. This ends the proof of statement (b).\\

It remains to prove  the third   statement in Theorem 1.2. Let $\Delta$ be an unknotted  circle in $\S$ and  let  $\Gamma$  be a spatial graph in $\S$ which either does not intersect $\Delta$ or intersects $\Delta$ exactly at one vertex $v$. Now, let $\tilde{\Gamma}$ be the pre-image of $\Gamma$ in the cyclic $p-$fold cover branched along $\Delta$. We consider the two maps   $\alpha$ and $\beta$  defined on the graph skein module of $\S$  as follows:
$\alpha(\Gamma)=Y(\tilde {\Gamma})(A)$ and $\beta(\Gamma)=Y(\tilde {\Gamma }!)(A)$, where $\tilde {\Gamma}!$ is the mirror
image of $\tilde {\Gamma}$. \\
Since  $Y(\tilde {\Gamma})(A)=Y(\tilde
{\Gamma})(A^{-1})$ \cite{Ya1}, then arguments similar to the ones used in the previous paragraph should enable us to prove that both of  the polynomials $\alpha$ and $\beta$
satisfy the following relations modulo $p$,
$$\alpha(\Gamma _{+} )(A)\equiv A^{4p}\alpha(\Gamma_{0})(A)+A^{-4p}\alpha(\Gamma
_{\infty})(A)-d^{p}\alpha(\Gamma _{\times}),
$$
and
$$\beta(\Gamma_{+})(A) \equiv A^{-4p}\beta(\Gamma_{0})(A)+A^{4p}\beta(\Gamma_{\infty})(A)-d^{p}\beta(\Gamma_{\times})(A),$$ where $\Gamma_{+}$, $\Gamma_{0}$, $\Gamma_{\infty}$ and  $\Gamma_{\times}$ are respectively the four graph diagrams  which appear  in the first
Yamada skein relation, see figure 5.\\
 It is obvious that if    $A^{4p}=A^{-4p}$, then $\alpha$ and $\beta$ are defined using the same
skein relations. Hence, $\alpha \equiv \beta$ modulo
$p, A^{8p}-1$. Since any spatial graph with   $\Z_p$-symmetry or a vertex-fixing $\Z_p$-symmetry can be constructed as $\tilde \Gamma$ for some spatial graph $\Gamma$, then we conclude that the congruence in statement (c) holds. This ends the proof of Theorem 1.2. \fin\\

\end{document}